\documentclass[12pt]{amsart}

\newtheorem{theorem}{Theorem}[section]

\theoremstyle{definition}

\theoremstyle{remark}

\numberwithin{equation}{section}

\newcommand{\be}{\begin{equation}}
\newcommand{\ee}{\end{equation}}

\newcommand{\EK}{E(K)}
\newcommand{\EKK}{E'(K)}
\begin{document}

\title{The number of points on an elliptic curve with square $x$-coordinates}


\author{Yu Tsumura}
\address{Department of Mathematics, Purdue University
150 North University Street, West Lafayette, Indiana 47907-2067
}

\email{ytsumura@math.purdue.edu}
\thanks{}

\subjclass[2010]{Primary 14H52; Secondary 11G20}

\date{}

\dedicatory{}

\begin{abstract}
Let $K$ be a finite field.
We know that  half of the elements of $K^*$ are square.
So it is natural to ask how many of them appear as $x$-coordinates of points on an elliptic curve over $K$.
We consider  a specific class of elliptic curves over  finite fields and show that  half of the $x$-coordinates on an elliptic curve are square.

\end{abstract}

\maketitle


\section{Introduction.}
Let $K$ be a finite field of characteristic different from $2$.
Let $a, b\in K$ satisfy $b\neq0$ and $r=a^2-4b\neq0$.
We consider a class of elliptic curves over $K$ given by $E: y^2=x^3+ax^2+bx$.

We know that  half of the elements of $K^*$ are square.
So the natural question is how many points on an elliptic curve have square $x$-coordinates.
Since  half of the elements of $K^*$ are square, we expect that  half of  the points on $\EK$ have square $x$-coordinate.

We show that this is the case when the elliptic curve is $E: y^2=x^3+ax^2+bx$.

\section{Result}
The main result is the following.
\begin{theorem}
Let $\EK$ be an elliptic curve defined above over a finite field $K$.
Let $S=\{P\in \EK | P=(x,y)$ with $x$ being square or $P=\infty \}$, where $\infty$ is the identity of $\EK$.
Then we have
\begin{align*}
 \#S= \left\{
\begin{array}{ll}
 \#\EK/2 & \textrm{ if } $b$ \textrm{ is a square in } $K$  \\
            \#\EK/2+1 & \textrm{ if } $b$ \textrm{ is not a square in } $K$ .
\end{array}\right.
\end{align*}
\end{theorem}
\begin{proof}
Consider $E': Y^2=X^3-2aX^2+rX$.
Then there is an isogeny $\phi:E'\longrightarrow E$ of degree $2$ defined by
\begin{equation}\label{equ:phi}
\phi: (X,Y) \longmapsto \left(\frac{Y^2}{4X^2},\frac{ Y(r-X^2)}{8X^2}\right).
\end{equation}
(See page 70, Example 4.5 in \cite{Silverman}.)

Then by the isomorphism theorem we have, $\EKK/{\rm Ker}(\phi) \cong {\rm Im}(\phi)$.
Since $E$ and $E'$ are isogenous, it is easy to see that $\#\EK=\#\EKK$.
(See page 153, Exercise 5.4 (a) in \cite{Silverman}.)
Also we have ${\rm Ker}(\phi)=\{ \infty,(0,0) \}$.
Hence we have
\begin{equation}\label{equ:iso}
\#{\rm Im}(\phi)=\#\EK/2.
\end{equation}

Now we show that 
\begin{align*}
 S= \left\{
\begin{array}{ll}
 {\rm Im}(\phi) & \textrm{ if } b \textrm{ is a square in } K  \\
             {\rm Im}(\phi) \cup \{(0,0)\} & \textrm{ if } b \textrm{ is not a square in } K.
\end{array}\right.
\end{align*}
Let $R$ denote the right hand side.
It is easy to see that the preimage of $(0,0)$ by $\phi$ over the algebraic closure $\bar{K}$ is $\{(a\pm2\sqrt{b}, 0)\}$.
Hence the preimage is in $\EK$ if and only if $b$ is a square in $K$. So $(0,0)$ is in both $S$ and $R$.

Suppose $P\in S$.
If $P$ is $\infty$, then clearly $P\in {\rm Im}(\phi)$.
If $P=(x,y)\neq(0, 0)$ with $x$ square, then we need to solve
\[\frac{Y^2}{4X^2}=x,\frac{ Y(r-X^2)}{8X^2}=y\]
for $(X,Y)\in \EKK$.
By the first equation, using the relation $Y^2=X^3-2aX^2+rX$, we have $X^2-2(a+2x)X+r=0$.
Hence we get $X=a+2x\pm 2\sqrt{x^2+ax+b}$.
Now since $y^2=x(x^2+ax+b)$ and $x$ is a nonzero square, we see that $x^2+ax+b$ is also a square.
Hence $X\in K$.
By the second equation for $Y$, we solve $Y=8yX^2(r-X^2)^{-1}\in K$.
Here $r-X^2\neq 0$ since $y\neq0$.
Hence we get a point $(X, Y) \in \EKK$ such that $\phi(X,Y)=(x,y)$.
Therefore, we have $S\subset R$.

The converse follows by  (\ref{equ:phi}) and we have $R=S$.

Now we have
\begin{align*}
 \#S= \left\{
\begin{array}{ll}
 \#{\rm Im}(\phi) & \textrm{ if } b \textrm{ is a square in } K  \\
             \#{\rm Im}(\phi)+1 & \textrm{ if } b \textrm{ is not a square in } K .
\end{array}\right.
\end{align*}
By the equation \ref{equ:iso},
\begin{align*}
 \#S= \left\{
\begin{array}{ll}
 \#\EK/2 & \textrm{ if } b \textrm{ is a square in } K  \\
             \#\EK/2+1 & \textrm{ if } b \textrm{ is not a square in } K.
\end{array}\right.
\end{align*}
\end{proof}

\bibliographystyle{amsplain}
\bibliography{square-x}

\end{document}